\documentclass[11pt,reqno,draft]{amsart}

\usepackage{amssymb}
\usepackage{mathrsfs}
\usepackage{cite}
\usepackage{amsfonts}
\newcommand{\rmv}[1]{}

\def\inv{^{-1}}
\def\p{\varphi}

\def\pinv{{\p \inv}}

\def\e<{\leq _{E}}

\def\malce{\mathbin{\hbox{$\bigcirc$\rlap{\kern-8.3pt\raise0,50pt\hbox{$\mathtt{m}$}}}}}

\def\1sk{^{(1)}}

\def\to{\rightarrow}

\def\Thmname{Theorem}
\def\Propname{Proposition}
\def\Lemmaname{Lemma}
\def\Definitionname{Definition}
\newtheorem{Thm}{\Thmname}
\newtheorem{Prop}[Thm]{\Propname}
\newtheorem{Lemma}[Thm]{\Lemmaname}
{\theoremstyle{definition}
}
{\theoremstyle{remark}
}

\title[Burnside problem for matrix semigroups]{Yet another solution to the Burnside problem for matrix semigroups}

\author{Benjamin Steinberg}
\address{School of Mathematics and Statistics\\
Carleton University \\
1125 Colonel By Drive\\
Ottawa, Ontario  K1S 5B6 \\
Canada}
\thanks{The author was supported in part by NSERC}
\email{bsteinbg@math.carleton.ca}
\date{\today}


\begin{document}
\begin{abstract}
We use the kernel category to give a finiteness condition for semigroups.  As a consequence we provide yet another proof that finitely generated periodic semigroups of matrices are finite.
\end{abstract}
\maketitle

\section{Introduction}
Schur proved that every finitely generated periodic group of matrices over a field is finite~\cite{curtis}.  McNaughton and Zalcstein established the corresponding result for semigroups~\cite{McZalc}.  Since then, a number of proofs of this result have appeared cf.~\cite{straubburn,guralnick,okninski,deluca,BerstelReutenauer}.  Here we use  ``global semigroup theory.''

\section{A finiteness condition for monoids}
Let $\p\colon M\to N$ be a homomorphism of monoids.  Define the \emph{trace} of $\p$ at $(n_1,n_2)\in N\times N$ to be the quotient  $\mathrm{tr}_{\p}(n_1,n_2)$ of the submonoid \[\{m\in M\mid n_1\p(m)=n_1, \p(m)n_2=n_2\}\] by the congruence $\equiv$ given by $m\equiv m'$ if, for all $m_i\in \pinv(n_i)$ ($i=1,2$), one has $m_1mm_2=m_1m'm_2$.  Recall that a monoid is called \emph{locally finite} if all its finitely generated submonoids are finite.  In this section, we prove the following finiteness result.

\begin{Thm}\label{tracefiniteness}
Let $\p\colon M\to N$ be a homomorphism of monoids such that $N$ is locally finite and  $\mathrm{tr}_{\p}(n_1,n_2)$ is locally finite for all $n_1,n_2\in N$.  Then $M$ is locally finite.
\end{Thm}

To prove this theorem, we use the kernel category from~\cite{Kernel} (see also~\cite[Section 2.6]{qtheor}). In this paper, we perform composition in small categories diagramatically: $fg$ means do $f$ first and then $g$. If $C$ is a category, then $C(c_1,c_2)$ will denote the hom set of arrows from $c_1$ to $c_2$.  We use $C(c)$ as short hand for the endomorphism monoid $C(c,c)$.  A category is said to be \emph{locally finite} if all its finitely generated subcategories are finite. First we need a lemma from~\cite{idempotentstabilizers}; the proof is essentially that of Kleene's Theorem.

\begin{Lemma}[Le Sa{\"e}c, Pin, Weil]\label{locallylocallyfinite}
Let $C$ be a category, all of whose endomorphism monoids are locally
finite.  Then $C$ is locally finite.
\end{Lemma}
\begin{proof}
Let $\Gamma$ be a finite graph with vertex set $V$ and let
$\p\colon \Gamma^*\to C$ be a functor that is injective on vertices (where $\Gamma^*$ denotes the free category on $\Gamma$).
For each $X\subseteq V$ and $v,v'\in V$, define $\Gamma^X_{v,v'}$ to
be the set of all paths in $\Gamma$ from $v$ to $v'$ visiting only
vertices from $X$  outside of its initial and terminal vertices. Then $\p(\Gamma^*) =
\bigcup_{v,v'} \p(\Gamma^{V}_{v,v'})$ and so it suffices to show each
$\p(\Gamma^{V}_{v,v'})$ is finite.  We proceed by establishing each
$\p(\Gamma^X_{v,v'})$ is finite by induction on $|X|$. Since
$\Gamma^{\emptyset}_{v,v'}$ is just the set of edges from $v$ to
$v'$ and $\Gamma$ is finite, the case $|X|=0$ is handled. Suppose
$X\neq \emptyset$ and choose $x\in X$.  Plainly,
\begin{equation}\label{locallocaleq}
\p(\Gamma^X_{v,v'}) = \p(\Gamma^{X\setminus \{x\}}_{v,v'})\cup
\p(\Gamma^{X\setminus \{x\}}_{v,x})\left(\p(\Gamma^{X\setminus
\{x\}}_{x,x})\right)^*\p(\Gamma^{X\setminus \{x\}}_{x,v'})
\end{equation}
 and by induction all sets $\p(\Gamma^{X\setminus
\{x\}}_{q,q'})$ are finite. Since the endomorphism monoid $C(\p(x))$ is
locally finite, we conclude $\left(\p(\Gamma^{X\setminus
\{x\}}_{x,x})\right)^*$ is also finite.  Finiteness of
$\p(\Gamma^X_{v,v'})$ now follows from \eqref{locallocaleq}.
 \qed\end{proof}

Let $\p\colon M\to N$ be a homomorphism of monoids.
Following~\cite{Kernel}, we define a small category $K_{\p}$,
called the \emph{kernel category} of $\p$.  The object set  of $K_{\p}$ is
$N\times N$. The arrows of $K_{\p}$ are equivalence classes $[n_1,m,n_2]$ of triples
$(n_1,m,n_2)\in N\times M\times N$ where \[[n_1,m,n_2]\colon (n_1,\p(m)
n_2)\to (n_1\p(m),n_2)\] and two triples $(n_1,m,n_2),(n_1,m',n_2)$ are identified if  $n_1\p(m)=n_1\p(m')$, $\p(m)n_2=\p(m')n_2$ and
$m_1mm_2=m_1m'm_2$ for all $m_i\in \pinv(n_i)$ with $i=1,2$.   Composition is given by
\[[n_1,m,\p(m') n_2][n_1\p(m),m',n_2] = [n_1,mm',n_2];\] the
identity at $(n_1,n_2)$ is $[n_1,1,n_2]$.  One easily verifies that $K_{\p}$ is a category~\cite{Kernel}  (see also~\cite[Section 2.6]{qtheor}).  Moreover, the endomorphism monoid $K_{\p}(n_1,n_2)$ is isomorphic to $\mathrm{tr}_{\p}(n_1,n_2)$.

It turns out that $K_{\p}$ is finitely generated whenever the domain is finitely generated and the codomain is finite.

\begin{Prop}\label{finitegeneration}
Suppose $\p\colon M\to N$ is a homomorphism of monoids with $M$ finitely generated and $N$ finite.  Then $K_{\p}$ is finitely generated.
\end{Prop}
\begin{proof}
Let $X$ be a finite generating set for $M$.
Then $K_{\p}$ is generated by all arrows of the form $[n_1,x,n_2]$ with $x\in X$.  Indeed, if $m=x_1\cdots x_r$, then
\begin{align*}
[n_1,m,n_2] =& [n_1,x_1,\p(x_2\cdots x_r)n_2][n_1\p(x_1),x_2,\p(x_3\cdots x_r)n_2]\\ & \cdots[n_1\p(x_1\cdots x_{r-1}),x_r,n_2],
\end{align*}
as required.
\end{proof}

We can now prove Theorem~\ref{tracefiniteness}
\begin{proof}[Proof of Theorem~\ref{tracefiniteness}]
Let $X$ be a finite subset of $M$, put $M'=\langle X\rangle$ and $N'=\p(M')$.  Then $N'$ is finite.  On the other hand, for $n_1,n_2\in N'$, clearly $\mathrm{tr}_{\p|_{M'}}(n_1,n_2)$ is a quotient of a submonoid of $\mathrm{tr}_{\p}(n_1,n_2)$ and hence locally finite.  Thus without loss of generality we may assume that $M$ is finitely generated and $N$ is finite.

From Proposition~\ref{finitegeneration} we conclude $K_{\p}$ is finitely generated.  Because each endomorphism monoid of $K_{\p}$ is of the form $\mathrm{tr}_{\p}(n_1,n_2)$, whence locally finite by hypothesis, Lemma~\ref{locallylocallyfinite} yields $K_{\p}$ finite.  Define a map $\psi\colon M\to K_{\p}$ by $\psi(m) = [1,m,1]$.  As $1\in \pinv(1)$, the equality $\psi(m)=\psi(m')$ implies that $m=1m1=1m'1=m'$ and hence $\psi$ is injective.  Thus $M$ is finite.
\end{proof}

\section{The Burnside problem for matrix semigroups}
By a periodic semigroup, we mean a semigroup so that each cyclic subsemigroup is finite.
Like many proofs of the McNaughton-Zalcstein Theorem, we begin with the case of an irreducible representation. We are not really doing anything new here;  our proof roughly follows~\cite{guralnick}.

\begin{Prop}\label{irreduciblecase}
Let $S$ be a finitely generated irreducible periodic subsemigroup of $M_n(K)$ where $K$ is an algebraically closed field.  Then $M$ is finite.
\end{Prop}
\begin{proof}
By a well-known theorem of Burnside, there are no proper irreducible subalgebras of $M_n(K)$ and hence we can find elements $s_1,\ldots, s_{n^2}$ of $S$ forming a basis for $M_n(K)$.  The trace form $(A,B)\mapsto \mathrm{tr}(AB)$ is a non-degenerate bilinear form on $M_n(K)$; let $s_1^*,\ldots,s_{n^2}^*$ be the corresponding dual basis.   Denote by $F$ the prime field of $K$.  Suppose $X$ is a finite generating set for $S$ and consider the finite set $A$ of elements $\mathrm{tr}(xs_is_j)$ and $\mathrm{tr}(xs_i)$ with $x\in X\cup\{1\}$, $1\leq i\leq n^2$. Put $E=F(A)$.  Since $S$ is periodic, $\mathrm{tr}(s)$ is either $0$ or a sum of roots of unity for $s\in S$.  Consequently, $E$ is a finite extension of $F$.

First we show that each $s_i^*$, for $i=1,\ldots,n^2$, can be written as a linear combination over $E$ of $s_1,\ldots, s_{n^2}$.  Indeed, let $C=(c_{ij})$ be the matrix given by $s_i^* = \sum c_{ij}s_j$ and   $D=(d_{ij})$ where $d_{ij} = \mathrm{tr}(s_is_j)$.  Then $D\in M_n(E)$ and $CD=I$ since
\[\sum c_{ik}d_{kj}=\sum c_{ik}\mathrm{tr}(s_ks_j) = \mathrm{tr}\left(\sum c_{ik}s_ks_j\right) = \mathrm{tr}(s_i^*s_j) =\delta_{ij}.\]  Thus $C=D^{-1}\in M_n(E)$, as required.

Now, for all $s\in S$, we can write
\begin{equation}\label{useful}
s= \sum_i \mathrm{tr}(ss_i^*)s_i = \sum_{i,j}c_{ij}\mathrm{tr}(ss_j)s_i.
\end{equation}
We claim that each element of $S$ can be written as a linear combination over $E$ of $s_1,\ldots, s_{n^2}$.  The proof is by induction on length.  For $x\in X$,  the claim is immediate from the definition of $E$ and \eqref{useful}.  Suppose $s=xs'$ with $x\in X$ and that $s'=\sum a_is_i$ with the $a_i\in E$.  Then $s=xs' = x\sum a_is_i$ and so $\mathrm{tr}(ss_j) = \sum a_i\mathrm{tr}(xs_is_j)\in E$.  An application of \eqref{useful} proves the claim.  As a consequence of the claim, it follows $\mathrm{tr}(s)\in E$ for all $s\in S$.

Let $T=\{\mathrm{tr}(s)\mid s\in S\}$.  We aim to prove that $T$ is finite.  Assuming this is true, it follows from \eqref{useful} that $S$ has at most $|T|^{n^2}$ elements.
We have two cases.  Suppose first that $K$ has characteristic $p>0$.  Then $E$ is a finite field and so trivially $T\subseteq E$ is finite.  Next assume that $K$ has characteristic $0$.  Let $k$ be the subfield obtained by adjoining to $\mathbb Q$ all entries of the elements of $X$.  Then $k$ is finitely generated over $\mathbb Q$ and $S\subseteq M_n(k)$.  The possible eigenvalues of an element $s\in S$ are zero and roots of unity satisfying a degree $n$ polynomial over $k$ (namely, the characteristic polynomial of $s$).  But a finitely generated extension field of $\mathbb Q$ has only finitely many such roots of unity
cf.~\cite[Proof of Theorem (36.2)]{curtis}.  Since the trace is the sum of the eigenvalues, it follows that also in this case $T$ is finite.
\end{proof}

The novel part of our proof is how we handle the reduction to the irreducible case.
The following lemma is a variant of a result from~\cite{AMSV}.

\begin{Lemma}\label{lemma:decomp}
Let $K$ be a ring with unit and $M\subseteq M_n(K)$ be a monoid of
block upper triangular matrices \[\begin{pmatrix} S & \ast \\
                                                 0 &
T\end{pmatrix}\ \text{with $S\subseteq M_m(K)$ and $T\subseteq M_r(K)$}.\]   Let $\p$ be the projection to the diagonal block and set $N=\p(M)$. Then, for all $n_1,n_2\in N$ the monoid $\mathrm{tr}_{\p}(n_1,n_2)$ embeds in the additive group $M_{m,r}(K)$.
\end{Lemma}
\begin{proof}
Fix $n_1,n_2\in N$ and put $R=\{m\in M\mid n_1m=n_1, mn_2=n_2\}$. Suppose that $n_1=(X,Y)$ and $n_2=(U,V)$.
We define a homomorphism $\psi\colon R\to M_{m,r}(K)$ as follows.   Given
\[m=\begin{pmatrix} A & B \\
                                                 0 &
C\end{pmatrix}\in R,\] define $\psi(m) = XBV$.
Note that
\begin{equation}\label{stabilizing}
XA=X,\ YC=Y,\ AU=U,\ CV=V
\end{equation}
by definition of $R$.  Using this we compute
\begin{equation}\label{kquotient}
\begin{pmatrix} X & Z\\
                  0 & Y\\\end{pmatrix}\begin{pmatrix} A & B \\
                                                 0 &
C\end{pmatrix}\begin{pmatrix} U & W \\
                                                 0 &
V\end{pmatrix} = \begin{pmatrix} XU & XW+XBV+ZV \\
                                                 0 &
YV\end{pmatrix}.
\end{equation}
Thus $\psi(m)=XBV$ determines and is determined by the right hand side of \eqref{kquotient} for any given $Z,W$. Therefore, $\psi(m)=\psi(m')$ if and only if $m$ and $m'$ represent the same element of $\mathrm{tr}_{\p}(n_1,n_2)$.

It remains to verify that $\psi$ is a homomorphism to the additive group  $M_{m,r}(K)$.  It clearly
sends the identity matrix to $0$. Also
if
\begin{equation*}
a=\begin{pmatrix} A & B \\
                                                 0 &
C\end{pmatrix},\ b= \begin{pmatrix} A' & B' \\
                                                 0 &
C'\end{pmatrix}\in R,
\end{equation*}
then $\psi(a) + \psi(b) = XBV + XB'V$, where as
\[ab=\begin{pmatrix} AA' & AB'+BC'\\
0  & CC'\end{pmatrix}.\] So $\psi(ab) = X(AB'+BC')V = XAB'V+XBC'V
= XB'V+XBV$ since $a,b\in R$ (cf.~\eqref{stabilizing}).
\end{proof}

We are now in a position to complete our proof of the theorem of McNaughton and Zalcstein~\cite{McZalc}.

\begin{Thm}[McNaughton and Zalcstein]
Let $S$ be a finitely generated periodic semigroup of $n\times n$ matrices over a field $K$.  Then $S$ is finite.
\end{Thm}
\begin{proof}
Without loss of generality we may assume that $K$ is algebraically closed.  Also assume that $S$ contains the identity matrix, since if it does not we may adjoin it.
We proceed by induction on $n$.  If $S$ is irreducible, we are done by Proposition~\ref{irreduciblecase}.   Otherwise, we can write \[S=\begin{pmatrix} M_1 & \ast \\ 0 & M_2\end{pmatrix}\text{where $M_1\subseteq M_m(K)$ and $M_2\subseteq M_r(K)$}\] are finitely generated periodic semigroups of matrices of strictly smaller degrees $m,r<n$.  By induction $M_1\times M_2$ is finite.  Consider the projection $\p\colon S\to M_1\times M_2$. Lemma~\ref{lemma:decomp} yields that $\mathrm{tr}_{\p}(x,y)$ is a periodic subsemigroup of the additive group $M_{m,r}(K)$, for any $x,y\in M_1\times M_2$.  In the case $K$ has characteristic $0$, this implies each such trace is trivial; if the characteristic of $K$ is $p>0$, then each such trace is an elementary abelian $p$-group.  In either case, it follows that each trace is locally finite and so Theorem~\ref{tracefiniteness} implies that $S$ is finite.
\end{proof}

\bibliographystyle{abbrv}
\bibliography{standard2}
\end{document}